\documentclass{proc-l}

\newtheorem{thm}{Theorem}[section]

\theoremstyle{definition}

\theoremstyle{remark}
\newtheorem{remark}[thm]{Remark}
\newtheorem{ex}[thm]{Example}

\numberwithin{equation}{section}

\newcommand{\bbr}{\begin{remark}}       
\newcommand{\eer}{\end{remark}}

\font\bbb=msbm10 scaled 1100

\newcommand{\pf}{{\em Proof: }}

\newcommand{\real}{\mbox{\bbb R}}
\newcommand{\comp}{\mbox{\bbb C}}

\newcommand{\zed}{\mbox{\bbb Z}} 
\newcommand{\eg}{{\em e.g.}}
\newcommand{\ie}{{\em i.e.}}

\newcommand{\ra}{\rightarrow}
\newcommand{\del}{\partial}
    
\newcommand{\tobyto}[4]                                         
        {\left[\begin{array}{cc}{#1}&{#2}\\{#3}&{#4}\end{array}\right]}

\begin{document}

\title{Tight contact structures via dynamics}
\author{John Etnyre}
\address{Department of Mathematics, Stanford University, 
	Palo Alto, CA, 94305}
\thanks{Supported in part by NSF Grant \# DMS-9705949.}

\author{Robert Ghrist}
\address{School of Mathematics, Georgia Institute of Technology,
	Atlanta, GA 30332}
\thanks{Supported in part by NSF Grant \# DMS-9508846.}


\begin{abstract}
We consider the problem of realizing tight contact structures on
closed orientable three-manifolds. By applying the theorems of Hofer
{\em et al.}, one may deduce tightness from
dynamical properties of (Reeb) flows transverse to the contact structure.
We detail how two classical constructions, Dehn surgery and branched 
covering, may be performed on dynamically-constrained links 
in such a way as to preserve a transverse tight contact structure.
\end{abstract}

\maketitle

\section{Contact geometry and dynamics}

For a more thorough treatment of the basic definitions and theorems 
related to the geometry and dynamics of contact structures  see, \eg, 
\cite{Aeb94}.


A {\em contact structure} $\xi$ on a 3-manifold $M$ is a totally 
non-integrable
2-plane field in $TM$.  More specifically, at each point $p\in M$ we have a
2-plane $\xi_p\subset T_pM$ that varies smoothly with $p$, with the property 
that $\xi$ is nowhere integrable in the sense of Frobenius: \ie, there exists 
(locally) a defining 1-form $\alpha$ (whose kernel is $\xi$) such 
that $\alpha\wedge d\alpha\not=0$. If $\alpha$ is globally defined, $\xi$ is 
called {\em orientable} and $\alpha$ a {\em contact 1-form} for $\xi$.  
We adopt the common restriction to orientable contact structures.

The interesting (and difficult) problems in contact geometry are all of 
a global nature: Darboux's Theorem (see, \eg, 
\cite{MS95,Aeb94}) implies that all contact structures are locally 
{\em contactomorphic}, or diffeomorphic preserving the plane fields. 
A similar result holds for a surface $\Sigma$ in a contact manifold $(M,\xi)$
as follows. Generically, $T_p\Sigma\cap \xi_p$ will be a line in $T_p\Sigma.$
This line field integrates to a singular foliation $\Sigma_\xi$ called the
{\em characteristic foliation} of $\Sigma$. One can show, as in the 
single-point case of Darboux's Theorem, that $\Sigma_\xi$ determines the germ 
of $\xi$ along $\Sigma$.

There has recently emerged a fundamental dichotomy in three dimensional
contact geometry.
A contact structure $\xi$ is {\em overtwisted} if there exists an
embedded disk $D$ in $M$ whose characteristic foliation $D_\xi$ contains
a limit cycle.  If $\xi$ is not overtwisted then it is called
{\em tight}.   Eliashberg \cite{Eli89} has completely
classified overtwisted contact structures on closed 3-manifolds --- the
geometry of overtwisted contact structures reduces to the algebra of
homotopy classes of plane fields.  Such insight into tight contact
structures is slow in coming. The only general method for constructing
tight structures is by Stein fillings (see \cite{Gom96,Eli90}) and the
uniqueness question has only been answered on $S^3$ \cite{Eli92}, $T^3$
\cite{Gir94,Kan95}, most $T^2$-bundles over $S^1$ \cite{Gir94}, 
and certain lens spaces $L(p,q)$ \cite{Etn97:lens}.

Thus we have the fundamental open question: {\em does every 3-manifold $M$
admit a tight contact structure?} Martinet \cite{Mar71} and Thurston and
Winkelnkemper \cite{TW75} have used surgery techniques to show that
all closed 3-manifolds admit contact structures.  However, 
their constructions do not yield tight contact structures.  
The current state of affairs is to be found in a recent theorem
of Eliashberg and Thurston \cite{ET97}, who show how to perturb 
taut foliations into tight contact structures. Then, the theorems of
Gabai \cite{Gab83} imply that any closed orientable 
irreducible 3-manifold with nonzero second Betti number 
$\beta_2$ supports a tight contact structure.

The outline of this paper is as follows: the remainder of this section 
consists of dynamical preliminaries and a recollection of the striking 
work of Hofer {\em et al.} concerning Reeb fields. We proceed in \S 2  
to carefully modify the well-known construction of Dehn surgery to 
preserve a tight contact structure when surgering over  
certain links. In \S 3, we turn to the 
procedure of branched covering and again show how to perform this
construction in such a way as to preserve a tight contact structure. 
In both cases, the link that the surgery / branching is performed on 
is constrained by the associated dynamics; thus, unfortunately, only 
certain tight manifolds are obtained by our methods. In particular, 
we do not surpass the existing theorems of Eliashberg and Thurston.
However, this marks the first examples of proving tightness of 
surgered contact structures without means of Stein filling. It is 
unknown whether the structures we construct are Stein fillable in 
general.

\subsection{The dynamics of Reeb flows}

A contact 1-form $\alpha$ carries
more geometry that does its contact structure $\xi=\ker\alpha.$  
In particular, given a contact form $\alpha$ there is a vector field 
$X$ uniquely determined by
$\alpha(X) = 1$ and $d\alpha(X,\cdot) = 0$.
The vector field $X$ is called the {\em Reeb vector field} \cite{Ree52}, and 
it encapsulates the ``extra geometry'' $\alpha$ carries, since the Reeb 
field is characterized by the properties of being transverse to $\xi$ and 
preserving the 1-form $\alpha$. 
In his recent work on the Weinstein conjecture \cite{Hof93} Hofer 
has found deep connections between the dynamics of the Reeb 
vector field $X$ and the tightness of $\xi$:
\begin{thm}[Hofer \cite{Hof93}]\label{thm_H}
Let $\xi$ be an overtwisted contact structures on the closed 3-manifold $M$.
Then the flow of the Reeb vector field associated to any contact 1-form
generating $\xi$ has at least one closed orbit of finite order in $\pi_1(M).$
\end{thm}
This can be refined by considering the dynamics of the closed orbits.
Following the standard usage \cite{HZ94}, a periodic orbit in a 
Hamiltonian flow is either {\em degenerate} or {\em nondegenerate}, 
depending on whether the spectrum associated to the linearized return 
maps for the orbits contains, or excludes respectively, one. The 
nondegenerate periodic orbits are either {\em elliptic} or 
{\em hyperbolic}, depending on whether these eigenvalues are on the unit 
circle or not respectively. 
\begin{thm}[Hofer, Wyzocki, and Zehnder \cite{HK97}]\label{thm_HWZ}
Let $\xi$ be an overtwisted contact structure on the closed 3-manifold $M$.
Then if the flow of the Reeb vector field associated to a contact 1-form
generating $\xi$ has no degenerate periodic orbits, then there exists 
at least one closed hyperbolic orbit of finite order in $\pi_1(M).$
\end{thm}
The proofs of the above theorems are highly nontrivial, relying 
primarily on Gromov's theory of pseudoholomorphic curves \cite{Gro85}.

\section{Dehn surgery on tight contact structures}\label{tsurg}

The operation of {\em Dehn surgery} is a very
efficient way of constructing closed orientable three manifolds.
A classical theorem in 3-manifold topology asserts that any closed
orientable three-manifold is obtainable via surgery on a link in 
$S^3$ \cite{Wal60,Lic62}. In this section, we show how
to preserve tightness under certain circumstances.

\subsection{Dehn surgery}

The object of Dehn surgery on a three-manifold is to drill out a
solid torus, and replace this with another solid torus inserted with
``twists.''

Let $\gamma$ denote a simple closed curve in $M^3$ having tubular
neighborhood $N$ diffeomorphic to $D^2\times S^1$, with
$\gamma$ as $\{0\}\times S^1$. Choose cylindrical coordinates
$(\rho,\theta,\phi)$ on $N$ such that the boundary curves
$m=\{(1,\theta,0)\}$ and $\ell=\{(1,0,\phi)\}$ correspond
respectively  to a meridian and a longitude of $N$. The meridional
curve is canonically defined, whereas the choice of the longitudinal
curve depends on the framing of the coordinate system.
Denote by $\Psi:\del N\ra\del N$ the diffeomorphism 
\begin{equation}\label{dehnmap}
        \Psi\left(\begin{array}{c}\theta \\ \phi\end{array}\right)
        = \tobyto{p}{s}{q}{t}
        \left(\begin{array}{c}\theta \\ \phi\end{array}\right) ,
\end{equation}
where $p, q, s,$ and $t$ are integers satisfying $pt-qs=1$.
The {\em $p/q$ Dehn surgery} of $M$ along $\gamma$ is performed
by removing $N$ from $M$ and regluing it via $\Psi$,
\begin{equation}
        M_{\gamma}(p/q) := \overline{M\setminus N}\bigcup_{\Psi} N ,
\end{equation}
resulting in the new manifold $M_{\gamma}(p/q)$ (completely 
determined by $\gamma, p$ and $q$).

\subsection{Model Contact Structures on $S^1\times D^2$}\label{model}

When performing Dehn surgery we will need to keep track of the Reeb
vector field in order to use Hofer's theorem to conclude our surgered 
manifold is tight. This is done by constructing model contact forms 
on $S^1\times D^2$. To this end, for any fixed $r>0$ we define the 
Hamiltonian function
$  H({\bf x}):=(x_1^2+y_1^2)+\frac{1}{r^2}(x_2^2+y_2^2) $
on $\real^4,$ where ${\bf x}:=(x_1,y_1,x_2,y_2)$.  The 1-form
$  \alpha:=\frac{1}{2}(x_1\, dy_1-y_1\, dx_1+x_2\, dy_2-y_2\, dx_2)$  
restricts to a tight contact 1-form $\hat{\alpha}$ on $\hat{S}=H^{-1}(1)$
\cite{Ben82}.   
We set $S:=\{(x_1,y_1,x_2,y_2) \vert x_1^2+y_1^2+x_2^2+y_2^2=1\}$ 
and define the map
$  \psi:S\ra \hat{S}:{\bf x}\mapsto {\bf x}/\sqrt{H({\bf x})}$.
Thus we obtain the tight contact structure $\alpha_H:=\psi^*(\hat\alpha)$ on
$S$. One may easily check that $\alpha_H=\frac{1}{H(x)}\alpha.$

Choose coordinates $0\leq\phi<2\pi$ on $S^1$ and polar coordinates
$(\rho,\theta)$ on $D^2.$  In these coordinates we define a map
\begin{equation}
  f:S^1\times D^2\ra S:(\phi,(\rho,\theta))\mapsto (\sin\rho\, e^{i\theta},
  \cos\rho\, e^{i\phi}).
\end{equation}
We define our distinguished contact forms on $S^1\times D^2$
as $\alpha_r:=f^*(\alpha_H),$ which in the above coordinates is
\begin{equation}\label{ctstrs}
  \alpha_r=\frac{1}{(\sin^2\!\rho+\frac{1}{r^2}\cos^2\!\rho)}
        (\sin^2\!\rho\, d\theta+\cos^2\!\rho\, d\phi).
\end{equation}
It is now simple to check
\begin{enumerate}
  \item For a fixed $\rho$ we get a torus $T_\rho$ in $S^1\times D^2$ by
	letting $\theta$ and $\phi$ vary.  The characteristic 
	foliation on $T_\rho$ is by lines with slope $-\tan^2\rho$.
  \item The Reeb vector field of $\alpha_r$ is tangent to the tori $T_\rho$
        	and has slope ${\frac{1}{r^2}}$ independent of $\rho.$
\end{enumerate}

\subsection{Tight surgery}

It is possible to perform Dehn surgery on a contact manifold
and obtain a new contact manifold \cite{Mar71}; 
however, without great care the contact
structure thus constructed will  be overtwisted. Using Stein
fillings, Eliashberg \cite{Eli90} and Gompf \cite{Gom96}
have shown how to build tight contact structures by {\em certain} 
surgeries on {\em any} knot in Stein fillable 3-manifolds (actually
symplectically semi-fillable would suffice, cf. 
\cite{Wei91}).  In this
section we show how to obtain a tight manifold by performing 
{\em any} surgery on a {\em certain} knot (\eg, the unknot in $S^3$).
\bbr\label{otremark}
To illustrate the difficulties in this type of construction consider the
following situation originally described to the authors by Maker-Limanov.
Let $\gamma$ be a closed transversal curve in $S^3$ (equipped
with its standard contact structure).  We can find a neighborhood
of $\gamma$ that is contactomorphic to $S^1\times D^2$ with the contact 
structure ker($d\phi+\rho^2 d\theta$), where$\phi$ is the 
angle coordinate on $S^1$ and $(\rho,\theta)$ are polar coordinates on $D^2$.  
Suppose that this neighborhood is large enough to include the torus $T$ 
formed by setting $\rho=\sqrt{2}.$  Now perform $-2/1$
Dehn surgery on $\gamma$ by cutting and pasting a solid torus not intersecting
$T.$  The characteristic foliation on $T$ is by $(-2,1)$ curves, thus they
bound disks in the surgered manifold.  Since the characteristic foliation
of $T$ has no singularities, it is not hard to find a disk with one of these
$(-2,1)$ curves as a limit cycle on its boundary. Thus the contact structure
one obtains is overtwisted.  When surgering tight contact
structures one must be careful to perform surgery on sufficiently 
``large'' tori.
\eer

\begin{thm}\label{tight_lens}
On $S^3$ with the (unique) tight contact structure, there exist 
tight ($p,q$)-Dehn surgeries on the unknot for arbitrary $p, q$.
\end{thm}
\pf 
For some irrational $r>0$, consider the contact form $\alpha_H$ 
on $S^3$ defined in \S 1.2. The Reeb field associated to 
$\alpha_H$ has precisely two periodic orbits, both of which are
elliptic. On a neighborhood of each of these closed orbits, the
contact form appears as in Equation~(\ref{ctstrs}); hence, 
the Reeb flow lies on invariant tori.
We will remove an invariant neighborhood of one of the 
periodic orbits, $\gamma$, 
and glue in a solid torus using $\Psi$ (from Equation~(\ref{dehnmap})),
thus performing a $p/q$ Dehn surgery on the unknot $\gamma$.  

Place coordinates $(\phi,(\rho,\theta))$ on $S^{1}\times D^{2}$ as in
Section~\ref{model}.  Recall that for fixed $\rho$ the torus
$T_{\rho}$ in $S^{1}\times D^{2}$ has as its characteristic foliation 
lines of slope $m_\rho=-\tan^2\rho$, 
and the Reeb vector field is tangent to $T_{\rho}$ with slope
$\frac{1}{r^{2}}.$
Pulling $\alpha_{r}$ back by $\Psi$ we obtain a contact structure on
a neighborhood of the boundary in $S^{1}\times D^{2}$ with characteristic
foliation on each $T_{\rho}$ of slope
\begin{equation}
n_\rho=-\left(\frac{pm_\rho-q}{sm_\rho-t}\right)
\end{equation}
and Reeb vector field tangent to $T_{\rho}$ with slope
\begin{equation}
\overline{r}=-\left(\frac{p-qr^2}{s-tr^2}\right).
\end{equation}
It is easy to check that given $p$ and $q$ one can find $s,t,r$ and
$\rho\in[0,\frac{\pi}{2})$ such that $pt-qs=1,$ $\overline{r}>0$, and 
$n_\rho<0.$  Now let $N=S^1\times D^2$, where the
$\rho$ variable is restricted to lie in the interval 
$[0,\tan^{-1}\sqrt{(n_\rho)}],$
and let $\alpha_{{1}/{\sqrt{\overline{r}}}}$ be the model contact form
constructed in \S \ref{model}.
One can now construct a map $\Phi$ from a neighborhood of the boundary of
$(N,\alpha_{1/\sqrt{\overline{r}}})$ to a neighborhood of the boundary of
$(S^1\times D^2,\Psi^*\alpha_r)$ that preserves the contact form
(to arrange this, make the map the identity on the invariant tori and
reparametrize in the $\rho$ direction so that the characteristic 
foliations on tori and the direction of the Reeb vector field are preserved). 
One may then use $\Psi\circ\Phi$ to glue 
$(N,\alpha_{{1}/{\sqrt{\overline{r}}}})$ 
to $(S^3\setminus (S^1\times D^2),\alpha_r)$. This contact form has a Reeb
field with neither degenerate nor hyperbolic periodic orbits; hence it 
is tight by Theorem~\ref{thm_H}
\qed

\bbr
Recall the lens space $L(p,q)$ is obtained from $S^3$ by performing
$-{p}/{q}$ Dehn surgery on an unknot. 
It has been known for a long time that all lens spaces admit tight contact
structures. Our interest in this theorem is the novel way of 
proving that these contact structures are tight --- using dynamical
properties of the Reeb vector fields to detect subtle geometric information.
This is the first ``surgery'' construction of tight contact structures 
that does not rely on Stein filling.
\eer

\bbr
It is not hard to compute the Euler class of the contact structure $\xi$
constructed in this example: let $D$ be the 2-skeleton
of the natural CW-decomposition coming from the surgery (which follows
since we are surgering an unknot). Then, using as the
generator of $H^{2}(L(p,q);\zed)$ the cochain that evaluates to
one on $D$, the Euler class of $\xi$ is $e(\xi)=q+1$. 
This follows from the formula for the Euler class in
\cite{Etn97:lens}, given that the characteristic foliation on $D$
has exactly one singular point (as can be seen using the local models).
When $p$ is even, there is a refinement of the Euler
class defined in \cite{Gom96} which may be likewise computed 
(see \cite{Etn97:lens} for a precise statement).
We note that in every case considered, this Euler class can be realized  by
a tight contact structure that can be Stein filled --- it is unknown 
whether this is true in general.
\eer

\bbr
In the above proof we never specifically addressed the problem of surgering
on ``sufficiently large tori'' discussed in Remark~\ref{otremark}.  It is
an interesting exercise to see that if one chooses $s,t,r$ and $\rho$ so that
$\overline{r}>0$ and $n_\rho<0,$ then a sufficiently large torus is being
surgered.
\eer

\section{Branched covers and tight contact structures}

Another way of building all three-manifolds is via
branched covers over knots and links. In this section, we show how one may
perform tight branched coverings of 3-manifolds along closed orbits of 
a suitable Reeb flow. 

\subsection{Branched covers over links}

To branch over a knot, one removes
a neighborhood of the knot, takes an $n$-fold cover of the complement,
and then fills in the tube(s) in such a way that the cross-sectional
map in the meridional direction is the $m$-fold singular cover
of the disc $D\subset\comp$ given by $z\mapsto z^m$.
More specifically, let $\gamma:=\{\gamma_i\}_1^n$
denote an $n$-component link in $M^{3}$.
Denote by $N$ a tubular neighborhood of $\gamma$ and by
$E:=\overline{M\setminus N}$ the exterior of $N$. For any
subgroup $G<\pi_1(E)$ of finite index, there is
a well-defined compact cover $p:M_G\ra E$.
Denote by $\{T_i\}$ the collection of boundary components
of $M_G$, each diffeomorphic to a torus. The cover $p$ restricts
to $p_i:={p}\vert_{T_i}$ on each torus.

Each boundary component of $E$ is a torus which may be fitted
with a meridian in such a way that each $p_i$
lifts this to a meridian for $T_i$ via an $m_i$-fold cover.
We may then construct $\tilde M$ by filling in the
$T_{i}$'s with $S^{1}\times D^{2}$'s
sending $\{\hbox{pt.}\}\times \partial{D^{2}}$ to the meridian.
After choosing a longitude for all the tori, each $p_{i}$ can be 
represented as $p_{i}(\theta,\phi)=(m\theta+k\phi, l\phi)$,
where $\theta$ and $\phi$ are the meridional and longitudinal coordinates,
respectively, and $m=m_{i}, k$, and $l$ are integers. If we extend each
$p_i:T^2\ra T^2$ to a map $\overline{p_i}:D^2\times S^1\ra D^2\times S^1$
via $(\rho,\theta,\phi)\mapsto(\rho,m\theta+k\phi,l\phi)$, then
the branched cover of $M$ over $\gamma$ via $G$ is defined to
be $\tilde M$ with projection
\begin{equation}     \tilde{p} : \tilde M\ra M \;\left\{\begin{array}{l}
                p \mbox{ on } M_G \\
                \overline{p_i} \mbox{ on each } D^2\times S^1
        \end{array} \right .
\end{equation}
Note the above projection map is not a smooth map since the
$\overline{p_i}$ are not smooth at $\rho=0.$  One could also define
the $\overline{p_i}$'s so that they are smooth: 
$(\rho,\theta,\phi)\mapsto(\rho^2,m\theta+k\phi,l\phi)$.  In this case, 
however,
$d\overline{p_i}=0$ at $\rho=0.$  We will make use of both the smooth and 
non-smooth versions in the following section.

\subsection{Tight branching over elliptic orbits}

In \cite{Gon87}, Gonzalo demonstrates lifting contact structures via 
a branched covering, and in this way also constructs contact structures
on all closed orientable 3-manifolds. There is no indication of 
tightness of such structures. In 
general, taking the (unbranched) cover of a tight contact
manifold can yield overtwisted contact manifolds \cite{Gom96} ---
so much more so for branched coverings.

We begin by showing how one can branch over certain elliptic periodic
orbits in a Reeb field to obtain tight contact structures.  We 
say a periodic orbit $\gamma$ in the Reeb flow of a contact form $\alpha$
is {\em locally integrable at $\gamma$} 
if there exist a neighborhood $N$ of $\gamma$ and 
(smooth) coordinates $(\rho,\theta,\phi)$ such that the Reeb field
takes the form $a(\rho)\frac{\partial}{\partial\theta} + 
b(\rho)\frac{\partial}{\partial\phi}$. 
These are precisely the action-angle coordinates from an integrable 
two degree-of-freedom Hamiltonian system, restricted to an energy surface.

\begin{thm}\label{thm_elliptic-branch}
Let $\alpha$ be a contact 1-form on $(M,\xi)$ such that the associated
Reeb vector field $X$ either (1) supports no closed orbits of finite
order in $\pi_1(M)$; or (2) supports no degenerate orbits and no hyperbolic 
orbits of finite order in $\pi_1(M).$  Moreover, assume that $X$ admits
a link of locally integrable periodic orbits $\gamma$. Then, any branched 
cover $\tilde{p}:\tilde{M}\ra M$ over $\gamma$ has a tight contact 
structure $\xi^p$ which is the lift of $\xi$ outside of a neighborhood of 
$\gamma$.
\end{thm}
\pf We assume without loss of generality that $\gamma$ is a 
single-component link. In order to pull back
the form $\alpha$ to a smooth contact form on the branched cover, we 
need to ensure that the form can be made locally $\theta$-equivariant.

{\em Case 1:} If the Reeb field $X$ has degenerate periodic orbits
near $\gamma$, then we may perturb the contact form as follows. 
Near any transverse loop such as $\gamma$, there exist coordinates for
which the contact structure is the kernel of $d\phi+\rho^2d\theta$
\cite[Thm. 8.3]{Aeb94}. In these
coordinates, $\alpha$ is of the form $f(\rho,\theta,\phi)(d\phi 
+ \rho^2d\theta)$ for some positive function $f$. 
To remove the $\theta$-dependence of $\alpha$
near $\gamma$, Taylor-expand $f$ as 
$f=f_0(\rho,\phi)+\overline{f}(\rho,\theta,\phi)$,
where $\overline{f}$ is $O(\rho)$. Then, choose a bump function 
$\chi(\rho)$ with support on $N$ attaining the value $1$ on a very small
neighborhood of $\rho=0$, and consider the form 
$\beta := g(\rho,\theta,\phi)(d\phi + \rho^2d\theta)$, where 
$g:=(f_0 + (1-\chi)\overline{f})$.
Since $\overline{f}$ is $O(\rho)$, the fact
that $f>0$ implies that $g>0$, and, hence, that $\beta$ is contact.

The Reeb field $Y$ for $\beta$ may have a very different periodic orbit
structure from $X$. Since $\gamma$ is locally integrable, the orbits of $Y$ 
are bound by invariant tori outside of a very small neighborhood. 
Hence, every closed orbit of $Y$ near $\gamma$ is a multiple of $\gamma$
in $\pi_1(M)$. It follows from hypothesis that $\gamma$ is of infinite order in 
$\pi_1(M)$, so it suffices to show that this multiple is always nonzero.
To do this, note that the $\frac{\partial}{\partial\phi}$-component
of $Y$ is given by $(2g+\rho g_\rho)/2g^2$. It suffices to show that the
numerator is nonzero on $N$. The first term, $2g$, is strictly positive. 
The second term, $\rho g_\rho$, may be made small through choice of
neighborhood and $\chi$, and hence does not 
overpower the (nonzero) $2g$ term. Thus every periodic orbit of $Y$ 
is also infinite order in $\pi_1$.

We may now branch since $\beta$ has the local normal form 
$\beta=f_0(\rho,\phi)(\rho^2d\theta + d\phi)$.
Pulling $\beta$ back by the non-smooth covering map $\overline{p_i}$ 
yields the local form
\begin{equation}
\tilde{\beta} = m\,f_0(\rho,l\phi)\rho^2\,d\theta + (k\,f_0(\rho,l\phi)\rho^2
		+l\,f_0(\rho,l\phi))\,d\phi ,
\end{equation}
on the branched cover. This form clearly extends over $\rho=0$, since 
$\beta$ was a smooth form. Thus $\tilde{\beta}$ is a
well-defined 1-form on $\tilde{N}$ which is a contact form since 
$\tilde{\beta}\wedge d\tilde{\beta} = (2ml\,f_0^2)\rho\,d\rho\wedge 
d\theta\wedge d\phi$. Moreover, since finite-order closed orbits on 
$\tilde{M}$ must project to finite-order closed orbits
on $M$, the Reeb field $Y$ of $\tilde{\beta}$ on $\tilde{M}$  
satisfies the hypotheses of Theorem~\ref{thm_H} and the contact structure 
$\xi^{p}=\hbox{ker}(\tilde{\alpha})$ is tight.


{\em Case 2:} If, in contrast, there are no degenerate periodic orbits 
near $\gamma$, then we may not perturb $\alpha$ to induce such. 
However, using the action-angle coordinates, we have that 
$X=a(\rho)\frac{\partial}{\partial\theta}+b(\rho)\frac{\partial}{\partial\phi}$
.
Since there are no degenerate periodic orbits, $a$ and $b$ are 
irrationally-related constants. In these coordinates, $\alpha$ takes on the 
form 
$\alpha = f\, d\theta + g\, d\phi + h\, d\rho$, where all the coefficients 
$f, g,$ and $h$ are functions of all three coordinates. 

By the definition of $X$, one has that $af+bg=1$ and 
$a(f_\rho-h_\phi)+b(g_\rho-h_\phi)=0$.
By differentiating the former with respect to each variable, one can
derive the equations $f_\phi = -\frac{a}{b}f_\theta$, $g_\phi = 
-\frac{a}{b}g_\theta$, and $h_\phi = -\frac{a}{b}h_\theta$. These
define first-order PDEs on $N$, and, in particular, on the invariant 
tori which foliate $N$. It is clear that the solutions to the Reeb field 
are characteristics of the PDEs; however, these
are dense on the invariant tori since $X$ has no degenerate closed
orbits. Hence, $f$, $g$, and $h$ are constants on each torus 
$\rho=c$ and these are all functions of $\rho$.

We thus have $\alpha$ of the form 
$\alpha=f(\rho)\,d\theta + g(\rho)\,d\phi + h(\rho)d\rho$.
Pulling $\alpha$ back by the non-smooth covering map $\overline{p_i}$ 
yields the local form
\begin{equation}
\tilde{\alpha} = m\,f(\rho)\,d\theta + (k\,f(\rho)
		+l\,g(\rho))\,d\phi + h(\rho)\,d\rho,
\end{equation}
on the branched cover. This form clearly extends over $\rho=0$, since 
$\alpha$ is a smooth form. Thus $\tilde{\alpha}$ is a
well-defined 1-form on $\tilde{N}$ which is a contact form since 
$\tilde{\alpha}\wedge d\tilde{\alpha} = ml\,\alpha\wedge d\alpha$.
Hence $\tilde{p}^{*}\alpha$ extends to a contact form on $\tilde{M}$. 
Moreover, since hyperbolic closed orbits on 
$\tilde{M}$ must project to (infinite-order) hyperbolic closed orbits
on $M$, the Reeb field of $\tilde{\alpha}$ on $\tilde{M}$  
satisfies the hypotheses of Theorem~\ref{thm_HWZ} and the contact structure 
$\xi^{p}=\hbox{ker}(\tilde{\alpha})$ is tight.
\qed



\begin{ex}[lens spaces]\label{ex_Lens2}
Consider the lens space $L(p,q)$ with the contact form $\alpha$
as constructed in Theorem~\ref{tight_lens}. The Reeb
vector field is an integrable field with precisely two closed orbits, 
$\gamma_{1}$ and $\gamma_{2}$,
which form the cores of a genus one Heegaard decomposition.  
As these orbits are elliptic, we may apply 
case (2) of Theorem~\ref{thm_elliptic-branch}.
The covers of $L(p,q)$ branched over $\gamma_{1}\cup
\gamma_{2}$ are of the form $L(p,q')$: it is an instructive
exercise to determine $q'$ for $p,q$ and the branching data. It
would appear that we have found more tight contact structures on
$L(p,q')$; however, it can be demonstrated that these structures are
all contactomorphic to the one constructed in
Theorem~\ref{tight_lens} (compute the Euler classes and then appeal to
the classification in \cite{Etn97:lens}).
\end{ex}

\begin{ex}[the three-torus]
The contact form $\alpha = (\sin z dx + \cos z dy) + 
\frac{1}{2}(\sin x dy + \cos x dz)$ has as its Reeb field (up 
to a nonzero rescaling)
\begin{equation}
	X = 2\sin z \frac{\partial}{\partial x} + 
		(\sin x + 2\cos z)\frac{\partial}{\partial y} +
		\cos x \frac{\partial}{\partial z} .
\end{equation}
This flow arises in the study of steady inviscid fluid flows 
\cite{Dom+86}. As this is a level set of an integrable Hamiltonian flow, 
it is simple to check that this vector field on $T^3$ has
no contractible closed orbits. The elliptic integral curves 
$\{(\pi/2, y, 0) : y\in\real/\zed\}$ and 
$\{(-\pi/2, y, \pi) : y\in\real/\zed\}$ are each a generator of 
$H_1(T^3)$ in the standard basis, and are locally integrable orbits. 
By Theorem~\ref{thm_elliptic-branch} branching over these 
curves yields tight contact structures on surface bundles over $S^1$. 
\end{ex}

\subsection{Tight branching over hyperbolic orbits}

We now consider branching over hyperbolic orbits.  This is a little more 
delicate and we need to make stronger global assumptions on the flow.

\begin{thm}\label{thm_anosov-branch}
Let $\alpha$ be a contact 1-form on $M$ such that the associated
Reeb vector field $X$ generates a structurally stable flow having no 
finite-order closed orbits (e.g., an Anosov flow).  Let $\gamma$ be any 
link of periodic orbits in the flow of $X$ and $\tilde{M}$ any branched 
cover over $\gamma$. Then $\tilde{M}$ has a tight contact form which
is the lift of $\alpha$ (outside of an arbitrarily small neighborhood of
$\gamma$).
\end{thm}

\pf
Consider a neighborhood $N$ of a component $\gamma_i$ of $\gamma.$
Using the smooth branching map pull back $\alpha\vert_N$ to a 1-form
$\beta$ on the cover $\tilde{N}.$  This smooth form $\beta$ is a contact 
form off of $\rho=0$.  Now set $\overline{\alpha}:=\beta+ \epsilon u(\rho) 
\rho^2 d\theta$, where $u(\rho)$ is a bump function with support 
on $\tilde{N}$ attaining $1$ near $\rho=0$, and $\epsilon$ is a small 
constant. It is not hard to check that for small $\epsilon$ 
the form $\overline{\alpha}$ is a contact form on all of $N$. 
Note $\rho=0$ is still a periodic orbit
of the Reeb field $\overline{X}$ for $\overline{\alpha}$.

The perturbation to the contact form, and hence the Reeb field, is 
equivariant with respect to the branching map.  Thus, away from  
$\rho=0$, the Reeb field of $\overline{\alpha}$ pushes down to a 
perturbation of the Reeb field of $\alpha$.  Thus flow lines of
the Reeb field of $\overline{\alpha}$ are mapped to flow lines of 
the perturbed field down stairs.  Moreover, since the Reeb field 
downstairs is structurally stable the perturbed field also has
no contractible orbits, implying the same for the Reeb flow of
$\overline{\alpha}$. 
\qed

\begin{ex}[pseudo-Anosov Reeb fields]
Let $M$ be the unit tangent bundle of a surface $\Sigma$ having 
constant negative curvature. The geodesic flow on $M$ is Anosov 
\cite{Ano67}, and preserves a transverse contact structure \cite{Pla72}. Let
$\alpha$ denote the natural contact form for which the flow is Reeb.
We may apply Theorem~\ref{thm_anosov-branch} to $(M,\alpha)$ to conclude: 
arbitrary branched covers over closed geodesics yield tight contact 
manifolds. This construction gives many interesting manifolds. 
\end{ex}

\bbr
The dynamics on the branched covers are no longer Anosov but can 
be lifted so as to be pseudo-Anosov (see, e.g. \cite{FLP79}). 
This provides a curious set of examples in light of the recent 
work of Benoist {\em et al.}, who show in \cite{BFL92} a strong 
rigidity among manifolds which admit an Anosov Reeb field. Namely, a  
Reeb field which is Anosov with $C^\infty$ splitting must be either a 
geodesic flow on a surface of constant negative curvature, or a certain 
time-reparametrization of this flow, or the lifted flow on an 
{\em unbranched} covering space of the unit tangent bundle. 
Our construction shows that relaxing the Anosov condition to a 
pseudo-Anosov condition greatly enlarges the 
class of 3-manifolds which admit such contact-preserving flows.
This presents an interesting problem in itself: {\em Classify which 
closed 3-manifolds admit a pseudo-Anosov Reeb flow with 
$C^\infty$ splitting}.
\eer

\providecommand{\bysame}{\leavevmode\hbox to3em{\hrulefill}\thinspace}

\end{document}